\documentclass{amsart}
\usepackage{latexsym, amssymb}
\usepackage{amsthm}

\def\beq{\begin{equation}}
\def\eeq{\end{equation}}\usepackage{url}

\newtheorem{lem}{Lemma}
\newtheorem{thm}[lem]{Theorem}
\newtheorem{prp}[lem]{Proposition}

\theoremstyle{definition}
\newtheorem{defn}{Definition}
\newtheorem{example}{Example}
\newtheorem{rem}{Remark}

\def\ra{\rightarrow}

\def\beqa{\begin{eqnarray}}
\def\eeqa{\end{eqnarray}}
\def\beqa{\begin{eqnarray}}
\def\eeqa{\end{eqnarray}}

\def\Out{\mathrm{Out}}
\def\Mod{\mathrm{Mod}}
\def\Comm{\mathrm{Comm}}

\def\vol{\mathrm{vol}}

\def\Isom{\mathrm{Isom}}

\begin{document}
\title[Isometries of symmetric spaces and moduli spaces]{Smith theory, $L^2$ cohomology, isometries of locally symmetric manifolds
and moduli spaces of curves}

\author{Grigori Avramidi}
\address{Dept. of Mathematics\\
5734 S. University Avenue\\
Chicago, Illinois 60637}
\email[G.~Avramidi]{gavramid@math.uchicago.edu}




\begin{abstract}
We investigate periodic diffeomorphisms of non-compact
aspherical manifolds (and orbifolds) and describe a 
class of spaces that have no homotopically trivial
periodic diffeomorphisms. Prominent examples are moduli spaces of curves
and aspherical locally symmetric spaces with $\chi\not=0$.
In the irreducible locally symmetric case, we show that 
no complete metric has more symmetry than the locally symmetric metric.
In the moduli space case, we build on work of Farb and Weinberger 
and prove an analogue of Royden's theorem for complete finite volume metrics.

\end{abstract}
\maketitle

\section{Introduction}
\label{intro}

A classical theorem of Borel, restricted to the case of closed hyperbolic
manifolds $(M,h_{hyp})$ of dimension $>2,$ asserts that for any Riemannian
metric $g$ on $M,$ the isometry group $\Isom(M,g)$ is abstractly isomorphic to a subgroup
of $\Isom(M,h_{hyp}).$
Thus, in a sense, the hyperbolic metric is the most
symmetric of all metrics. In this paper, we study this and related
phenomena for some noncompact locally symmetric spaces and for
moduli spaces of curves (with any number of punctures). Our work gives new 
information, for instance, about the isometry group of the McMullen metric
on Teichm\"uller space and establishes a version of Royden's calculation.
It also extends Borel's result to noncompact locally symmetric spaces with $\chi\not=0$. 

A context for this work is suggested by the recent geometric approach to 
computing the isometry groups of closed aspherical manifolds and their universal 
covers developed by Farb and Weinberger in \cite{farbweinbergerisometries}.
In the case of a closed hyperbolic manifold, 
this approach gives a dichotomy for the isometry group
of the universal cover. Either $g$ is a constant multiple
of the hyperbolic metric, or the isometry group of the universal cover is 
discrete and contains the covering translations as a subgroup of index 
$\leq\vol(M,h_{hyp})/\varepsilon(n)$, where the constant $\varepsilon(n)$
is the volume of the smallest $n$-dimensional hyperbolic orbifold.

The object of this paper is to extend the simple picture
described above for closed hyperbolic manifolds to two classes of 
{\it non-compact} aspherical orbifolds: the locally symmetric spaces of non-zero
Euler characteristic, and moduli spaces of curves. 
A fundamental question in the context of this problem is whether 
there are nontrivial periodic diffeomorphisms of the universal cover
which commute with the action of the fundamental group.
Closely related is the following 
\subsection*{Question:} Is every homotopically trivial periodic diffeomorphism 
of the aspherical manifold $M$ equal to the identity?
\subsection{Locally symmetric spaces}
We show that a locally symmetric manifold $M$ of non-zero Euler characteristic
has no homotopically trivial periodic diffeomorphisms. 
Combining this with Nielsen realization for surfaces, or Margulis-Mostow-Prasad rigidity
in higher dimensions, we see that the locally symmetric metric is the most symmetric one.
That is, the isometry group $\Isom(M,g)$
of any other complete Riemannian metric $g$ on $M$ is isomorphic to a subgroup of $\Isom(M,h_{sym}).$
If $g$ has finite volume, then we also get a result for the isometry group of the universal cover.

\begin{thm}
\label{symmetric}
Let $(M,h_{sym})$ be a non-positively curved, finite volume, irreducible locally symmetric manifold.
Suppose that $\chi(M)\not=0.$ If $g$ is any complete finite volume Riemannian metric on $M$ then either
\begin{itemize}
\item
$g$ is a constant multiple of the locally symmetric metric, or
\item
the isometry group of the universal cover $\Isom(\widetilde M,\widetilde g)$ is discrete and contains
$\pi_1M$ as a subgroup of index $\leq\vol(M,h_{sym})/\varepsilon(\widetilde h_{sym}).$
\end{itemize}
\end{thm}
The constant $\varepsilon(\widetilde h_{sym})$ is the volume of the smallest locally symmetric orbifold 
covered by $(\widetilde M,\widetilde h_{sym}).$ It does not depend on the metric $g.$ 
The non-positively curved symmetric spaces whose lattices have non-zero Euler characteristics include even-dimensional
hyperbolic spaces, Hermitian symmetric spaces with no compact or Euclidean factors, and products of these. 

\subsection{Moduli spaces of surfaces}
Let $\mathcal M_{g,n}^{\pm}$ be the moduli space of hyperbolic structures
on a genus $g$ orientable 
surface 
with $n$ punctures\footnote{This is double covered by the space of oriented hyperbolic structures $\mathcal M_{g,n}$.},
$X$ a smooth manifold and $X\ra\mathcal M_{g,n}^{\pm}$ a finite degree regular branched cover. 
We show that $X$ cannot have any non-trivial, homotopically trivial, periodic
diffeomorphisms. 
If the real dimension of the moduli space $\mathcal M_{g,n}^{\pm}$ is at least six, 
then work of Ivanov\cite{ivanov} implies that every periodic diffeomorphism of 
$X$ is homotopy equivalent to a covering translation.  
Together, these facts imply that   
the moduli space $\mathcal M_{g,n}^{\pm}$ is a {\it minimal} orbifold.
It is not a finite branched cover of any smaller orbifold.

\begin{rem}
The four-dimensional moduli spaces are related by a degree five cover $\mathcal M_{1,2}^{\pm}\ra \mathcal M_{0,5}^{\pm}$
corresponding to the quotient of the twice-punctured torus by the hyperelliptic involution.
In this case, we show that the orbifold $\mathcal M_{0,5}^{\pm}$ is minimal.  
\end{rem}

Combining this with a result of Farb and Weinberger (1.2 in \cite{farbweinbergerroyden}), 
which shows for a complete finite volume
metric on moduli space that the isometry group of the universal cover 
contains the covering translations as a {\it finite index} subgroup, we establish 
an exact analogue of Royden's computation of isometries of the Teichm\"uller metric. 
\begin{thm}
\label{generalroyden2}
Suppose that $S$ is not a sphere with $\leq 4$ punctures or a torus with $\leq 2$ punctures.
Let $\widetilde g$ be a complete Finsler metric on Teichm\"uller space $\mathrm{Teich}(S)$ which is invariant
under the extended mapping class group $\Mod^{\pm}(S)$ and has finite covolume. 
Then, all isometries of the metric $\widetilde g$ are elements of the extended mapping class group
acting by covering translations.
\end{thm}
For the Teichm\"uller metric this result was originally
proved in \cite{royden} using analytic methods, but for McMullen's K\"ahler hyperbolic metric and 
other finite volume complete metrics on moduli space it is new. 

Curiously, an infinite covolume metric on Teichm\"uller space can have extra isometries. 
We give an example of a metric whose isometry group contains the free product of $\mathbb Z$ with a 
finite index subgroup of the mapping class group. More generally, we have the following 
\begin{example}
Given two non-compact aspherical 
manifolds $M$ and $N$ whose universal covers
are diffeomorphic to $\mathbb R^m$ there is a complete infinite volume Riemannian
metric $g$ on $M$ so that the free product $\pi_1M\star \pi_1N$ acts by isometries
on the universal cover $(\widetilde M,\widetilde g).$
\end{example}

\subsection{The general picture}
Theorems \ref{symmetric} and \ref{generalroyden2} are both special cases of a single method,
whose range of applicability we will now describe. 
An aspherical manifold can have homotopically trivial periodic diffeomorphisms if the fundamental group
has a center, or if it is `too small' compared to the dimension of the manifold. 
For instance, the line $\mathbb R$ has an involution, while the plane $\mathbb R^2$ and the circle $S^1$
both have rational rotations. 
It is known that aspherical manifolds of non-zero Euler characteristic have centerless fundamental groups 
(7.2 in \cite{luck}). To describe a class of manifolds with `large' fundamental groups, we need the following notion.
\begin{defn} 
A manifold $\widetilde M$ with a properly discontinuous $\Gamma$-action is {\it $\Gamma$-tameable} 
if there is complete Riemannian bounded geometry\footnote{The sectional curvatures are bounded above and below, and
the injectivity radius is $\geq \varepsilon>0$.} metric on $\widetilde M$
which is $\Gamma$-invariant and has finite $\Gamma$-covolume. Sometimes we will
call such a manifold {\it $\Gamma$-tame}. 
\end{defn}
Manifolds of this sort were studied by Cheeger and Gromov in \cite{cheegergromov,cheegergromovgroups} 
and \cite{cheegergromovbounds}. They behave like $\Gamma$-covers of closed manifolds 
from the point of view of $L^2$ cohomology.
Examples are regular covers of closed manifolds, symmetric spaces 
with the action of a lattice and---by McMullen's result \cite{mcmullen}---Teichm\"uller spaces
with the action of the mapping class group. 

The main technical result of this paper is the following theorem.
It describes conditions on an aspherical manifold which ensure there are 
no homotopically trivial periodic diffeomorphisms. 
\begin{thm}
\label{noactions2}
Let $M$ be an aspherical manifold with $\pi_1M$-tame universal cover.
Suppose all three of the following conditions hold.
\begin{enumerate}
\item
The fundamental group $\pi_1M$ is residually finite.
\item
$M$ is the interior of a compact manifold with boundary\footnote{This assumption is convenient, 
but not necessary. See subsection \ref{infinitetype}.} 
$\partial M.$
\item
The Euler characteristic $\chi(M)$ is non-zero.
\end{enumerate}
Then,
\begin{enumerate}
\item
Any homotopically trivial periodic diffeomorphism of $M$ is the identity.
\item
$\pi_1M$ does not commute with any non-trivial compact Lie group $K<\mathrm{Diff}\widetilde M.$
More generally, $\pi_1M$ does not normalize any such $K$. 
\end{enumerate}
\end{thm} 
Now pick a complete Riemannian (or more generally Finsler) metric $g$ on $M$.
The first application of Theorem \ref{noactions2} is a non-compact analogue
of Borel's result.
\begin{thm}
\label{maintheorem}
Under the hypotheses of Theorem \ref{noactions2}, the isometry group $\Isom(M,g)$ is isomorphic to a finite subgroup of $\Out(\pi_1M).$ 
\end{thm}

If the metric $g$ has finite volume, then we  
also get the following recognition theorem for irreducible $\chi\not=0$ locally symmetric spaces
among aspherical manifolds. 

\begin{thm}
\label{recognizinglocsym}
Suppose $M$ satisfies the hypotheses of Theorem \ref{noactions2}, and in addition, 
$\pi_1M$ is irreducible\footnote{That is, if $A\times B$ is a finite index subgroup of 
$\Gamma$ then either $A$ or $B$ is trivial.}. 
If $g$ is a finite volume complete Riemannian (or Finsler) metric on $M,$ then
\begin{itemize}
\item
$(M,g)$ is an irreducible, non-positively curved, finite volume, $\chi\not=0$ locally symmetric space,
or
\item
The isometry group $\Isom(\widetilde M,\widetilde g)$ is discrete and contains $\pi_1M$ as a finite index subgroup.
\end{itemize}
\end{thm}

For closed aspherical manifolds, our Theorem \ref{recognizinglocsym} is (the $\chi\not=0$ case of) 
Theorem $1.3$ of \cite{farbweinbergerisometries}. For moduli spaces, it is Theorem $1.2$ of 
\cite{farbweinbergerroyden}. For noncompact $\chi\not=0$ locally symmetric spaces, it is new.

The paper is organized as follows.
In section \ref{lifts} we review the covering space theory needed to lift isometries to the universal cover.
In section \ref{Smith theory} we recall what we need about fixed point sets of $\mathbb Z/p$ actions.
Then, we describe a strengthening of the basic Smith theory result which applies to our situation.
This is enough to give a proof of Borel's theorem in section \ref{Borel theorem}. 
In section \ref{intersections} we formulate precisely when a manifold has `many intersections'
(condition ($\pitchfork$)) and show that such a manifold can't be a bundle.
In section \ref{L2} we recall what we need about Betti numbers of covers
and construct finite covers satisfying ($\pitchfork$). We also explain how to remove
the hypothesis that $M$ is the interior of a compact manifold with boundary.
In section \ref{noactions} we prove Theorem \ref{noactions2}.
In section \ref{recognition} we prove Theorem \ref{recognizinglocsym}.
As corollaries of this, we prove Theorems \ref{generalroyden2} and \ref{symmetric}
in sections \ref{royden} and \ref{locsym}, respectively.
In section \ref{infinitevolume} we show that Theorem \ref{recognizinglocsym} 
does not extend to metrics of infinite covolume.
In section \ref{BorelStar} we prove Theorem \ref{maintheorem}.

\subsection*{Acknowledgements} I would like to thank my advisor Shmuel Weinberger
for his guidance, inspiration and encouragement,
Benson Farb for suggesting that a lack of homotopically trivial isometries should yield Royden's theorem,
and T$\hat{\mathrm{a}}$m Nguy$\tilde{\hat{\mathrm{e}}}$n Phan 
for helpful discussions about locally symmetric spaces. I would also like to thank all three
for reading earlier versions of the paper and making a number of suggestions. 

\section{\label{lifts}Lifting isometries}

Let $f:M\ra M$ be a homeomorphism of the manifold $M$.
A {\it lift} of this homeomorphism to the universal cover $\widetilde M\ra M$ is a homeomorphism
$\widetilde f$ making the diagram

\begin{equation}
\begin{array}{ccc}
\widetilde M&\stackrel{\widetilde f}\longrightarrow &\widetilde M\\
\downarrow& &\downarrow\\
M&\stackrel{f}\longrightarrow &M
\end{array}
\end{equation}
commute. Composing $\widetilde f$ with a covering translation $\gamma\in\pi_1M$ gives
another lift $\gamma\widetilde f$ of $f.$ If $G$ is a group acting by homeomorphisms on the manifold
$M$ then the lifts $L$ of these homeomorphisms to the universal cover form a group extension
\begin{equation*}
\label{liftsextension}
1\ra\pi_1M\ra L\ra G\ra 1.
\end{equation*}
It is well known (see 8.8. in \cite{maclane}) that any group extension 
\begin{equation*}
1\ra H\ra L\ra G\ra 1
\end{equation*}
is determined by
\begin{enumerate}
\item
a representation 
$\rho:G\ra\Out(H),$ and
\item
a class in the cohomology group\footnote{The coefficients $Z(H)$ are a $G$-module via $\rho$.} $H^2(G;Z(H))$. 
\end{enumerate}
In the special case when the homomorphism $\rho$ and the center $Z(H)$ are both trivial
we get the trivial extension. 
\begin{prp}
\label{productlift}
If $1\ra H\ra L\ra G\ra 1$ is a group extension, $H$ has trivial center and the conjugation
homomorphism $\rho:G\ra\Out(H)$
is trivial, then the extension $L$ is isomorphic to the product $H\times G.$
\end{prp}
\begin{proof}
Since the homomorphism $\rho$ is trivial, for every $g\in L$ there is $\gamma(g)$ in $H$ such
that conjugation by $g$ agrees with conjugation by $\gamma(g)$ on $H.$ The element $\gamma(g)$ is 
unique because $H$ has trivial center, so we get a homomorphism $L\ra H, g\mapsto \gamma(g)$
which splits the extension. 
\end{proof}


\section{\label{Smith theory}$\mathbb Z/p$-actions and Smith theory}\

Let $p$ be a prime.
In this section we describe the homological structure of the fixed point set of a $\mathbb Z/p$-action. 
The basic theme is that the fixed point set cannot be $\mathbb Z/p$-homologically 
more complicated than the ambient space. 
First, suppose that a contractible manifold $\widetilde M$ is equipped with a complete, non-positively curved metric $\widetilde g$ 
and the group $\mathbb Z/p$ acts by isometries of this metric. Then the fixed point
set $\widetilde F$ is a totally geodesic submanifold.
In particular, it is non-positively curved and has no closed geodesics, so it is also contractible.
For a general metric we can only conclude that the
fixed point set `looks contractible from the point of view of $\mathbb Z/p$-homology.'
More precisely, one has the following result of Smith theory.

\begin{thm}[$2.2$ in chapter VII of \cite{bredon}]
\label{Smith1}
Suppose the group $\mathbb Z/p$ acts smoothly on a smooth, contractible manifold $\widetilde M.$
Then, the fixed point set $\widetilde F\subset \widetilde M$ of the $\mathbb Z/p$-action
is a $\mathbb Z/p$-acyclic, smooth submanifold.\footnote{In particular, the fixed point set is connected.} That is,

\begin{equation}
H_*(\widetilde F;\mathbb Z/p)\cong H_*(\widetilde M;\mathbb Z/p).
\end{equation}
\end{thm}

\subsection{\label{Equivariant Smith}$\mathbb Z/p\times\pi_1M$-actions}
The $\mathbb Z/p$-actions on the universal cover $\widetilde M$ we are interested in are obtained
by lifting homotopically trivial $\mathbb Z/p$-actions on the base manifold $M$. These actions
{\it commute with the covering action of the fundamental group $\pi_1M$.} (See Proposition \ref{productlift}.)
In this situation, we have an amplification of Theorem \ref{Smith1}.

\subsection*{Notation} For the rest of this section, $M$ is a smooth aspherical manifold 
and $\widetilde M$ is its universal cover. We are given a $\mathbb Z/p$-action on the universal cover 
which commutes with the $\pi_1M$-action, hence descends to a $\mathbb Z/p$-action on the manifold $M$.
We pick a triangulation of the smooth manifold $M$ in which 
the smooth $\mathbb Z/p$-action is simplicial (this can be done by 7.1 in \cite{illman}.) 
The universal cover $\widetilde M$ is given the lifted (hence $\pi_1M$-equivariant) triangulation. 
All chain complexes are simplicial. 

\subsection{Projecting the $\mathbb Z/p$-action to covers}
Let $\Gamma<\pi_1M$ be a subgroup of the fundamental group and and let $\phi\in\mathbb Z/p.$
Since the $\mathbb Z/p$-action on the universal cover $\widetilde M$ commutes with the $\pi_1M$-action, it 
descends to any other cover $M'=\widetilde M/\Gamma$ via the formula

\begin{equation}
\phi(\Gamma x)=\Gamma\phi(x).
\end{equation}
The following lemma shows that the fixed point set $F'$ of the $\mathbb Z/p$-action on $M'$ is the projection
$\widetilde F/\Gamma$ of the fixed point set $\widetilde F$ of the $\mathbb Z/p$-action on the universal cover $\widetilde M.$ 
In other words, we have the commutative diagram 

\begin{equation}
\label{fixedsetsandcovers}
\begin{array}{ccc}
\widetilde F&\hookrightarrow&\widetilde{M}\\
\downarrow& &\downarrow\\
F'&\hookrightarrow&M'.
\end{array}
\end{equation}
The vertical maps are $\Gamma$-covers and the horizontal maps are inclusions.

\begin{lem}
If $\phi$ preserves the $\Gamma$-orbit of a point $x\in\widetilde M,$ then $\phi$ fixes the orbit pointwise.
\end{lem} 

\begin{proof}
If $\phi(\Gamma x)=\Gamma x,$ then there is 
$\gamma\in\Gamma$ such that 
\begin{equation}
\phi(x)=\gamma x.
\end{equation} 
Since $\phi$ and $\gamma$ commute and $\phi\in\mathbb Z/p$, we find

\begin{equation}
x=\phi^p(x)=\phi^{p-1}\phi(x)=\phi^{p-1}\gamma x=\gamma\phi^{p-1}(x)=\cdots=\gamma^px.
\end{equation}
Thus, the element $\gamma^p\in\Gamma\subset\pi_1(M)$ fixes the point $x\in\widetilde M.$ 
Since the fundamental group acts freely on $\widetilde M$, we conclude $\gamma^p=1.$
The fundamental group of the aspherical manifold $M$ is torsionfree, so we must have $\gamma=1.$
Consequently, $\phi(x)=x$ and for any $\tau\in\Gamma$

\begin{equation}
\phi(\tau x)=\tau\phi(x)=\tau x,
\end{equation}  
so $\phi$ fixes the entire orbit $\Gamma x$ pointwise. 
\end{proof}

\subsection{Homology and cohomology with coefficients in a $\mathbb Z/p[\pi_1M]$-module}
Next, we will see that the inclusion of the fixed point
set $F'\subset M'$ is a $\mathbb Z/p$-homology equivalence in every cover $M'\ra M.$

\begin{thm}
\label{coeffsmith}
Let $M$ be a smooth aspherical manifold and $\widetilde M$ its universal cover.
Suppose we have a smooth $\mathbb Z/p$-action on $\widetilde M$ that commutes with the action
of $\pi_1M$ by covering translations. Let $F\subset M$ be the fixed point set of the projected $\mathbb Z/p$-action
on $M$. Let $V$ be a $\mathbb Z/p[\pi_1M]$-module. Then, the inclusion of the fixed point set induces
isomorphisms on homology and cohomology with coefficients in $V$, i.e.

\begin{eqnarray}
H_*(F;V)&\cong &H_*(M;V),\\
H^*(M;V)&\cong &H^*(F;V).
\end{eqnarray}
\end{thm}

\begin{proof}
The key point is that 

\begin{equation}
\label{exact}
0\ra C_*(\widetilde F;\mathbb Z/p)\ra C_*(\widetilde M;\mathbb Z/p)\ra C_*(\widetilde M,\widetilde F;\mathbb Z/p)\ra 0
\end{equation}
is an exact sequence of complexes of free $\mathbb Z/p[\pi_1M]$-modules, 
and the complex $C_*(\widetilde M,\widetilde F;\mathbb Z/p)$ is acyclic.
\begin{itemize}
\item
The sequence is exact by definition.
\item
It is a sequence of complexes of $\mathbb Z/p[\pi_1M]$-modules because the $\mathbb Z/p$-action commutes with the $\pi_1M$-action.
\item
For each $k$-simplex $\sigma_k$ in the base $M$, pick a lift $\widetilde\sigma_k$ in the universal cover $\widetilde M.$
Then, 
\begin{eqnarray*}
C_k(\widetilde F;\mathbb Z/p)&\cong&\oplus_{\sigma_k\in F}\hspace{0.6cm}\mathbb Z/p[\pi_1M]\widetilde\sigma_k,\\
C_k(\widetilde M;\mathbb Z/p)&\cong&\oplus_{\sigma_k}\hspace{1cm}\mathbb Z/p[\pi_1M]\widetilde\sigma_k,\\
C_k(\widetilde M,\widetilde F;\mathbb Z/p)&\cong&\oplus_{\sigma_k\notin F}\hspace{0.6cm}\mathbb Z/p[\pi_1M]\widetilde\sigma_k,
\end{eqnarray*}
show that the $\mathbb Z/p[\pi_1M]$-modules in the sequence (\ref{exact}) are all free.
\item
Theorem \ref{Smith1} shows that the complex $C_*(\widetilde M,\widetilde F;\mathbb Z/p)$ is acyclic. 
(Its homology vanishes in all dimensions.)
\end{itemize}
Recall that homology and cohomology with coefficients in $V$ are computed by\footnote{Here  $\widetilde{}$  denotes
the $\pi_1M$-cover.}

\begin{eqnarray}
H_*(-;V)&=&H_*(C_*(\widetilde-;\mathbb Z/p)\otimes_{\mathbb Z/p[\pi_1M]}V),\\
H^*(-;V)&=&H_*(\mathrm{Hom}_{\mathbb Z/p[\pi_1M]}(C_*(\widetilde-;\mathbb Z/p),V)).
\end{eqnarray} 
Since the $\mathbb Z/p[\pi_1M]$-modules
appearing in the sequence (\ref{exact}) are all free, the 
functor $-\otimes_{\mathbb Z/p[\pi_1M]}V$ preserves exactness of the sequence (\ref{exact}) and acyclicity
of the relative complex $C_*(\widetilde M,\widetilde F;\mathbb Z/p),$ so the long exact homology sequence
shows that $H_*(F;V)\ra H_*(M;V)$ is an isomorphism. The same remarks for the functor $\mathrm{Hom}_{\mathbb Z/p[\pi_1M]}(-,V)$
give the cohomology isomorphism.
\end{proof}
Note, in particular, that if $\Gamma<\pi_1M$ is a subgroup then homology with coefficients in the module
$V=\mathbb Z/p[\pi_1M/\Gamma]$ is just homology in the cover $M'=\widetilde M/\Gamma$, so the inclusion
of the fixed point set is a $\mathbb Z/p$-homology isomorphism

\begin{equation}
H_*(F';\mathbb Z/p)\ra H_*(M';\mathbb Z/p)
\end{equation}
in every cover $M'\ra M.$ This case is also proved by Conner and Raymond in the appendix of \cite{connerraymond}.

\subsection*{$\mathbb Z/p^N$-equivalences.}
Later on, we'll need the standard fact that a $\mathbb Z/p$-homology equivalence is also a $\mathbb Z/p^N$-homology
equivalence for any $N\geq 1.$ 
For a map $f:X\ra Y$, let $C_f$ be its mapping cone.
The map $f:X\ra Y$ is an $R$-homology equivalence if and only if its mapping cone is $R$-acyclic, i.e. $\tilde H_*(C_f;R)=0.$

\begin{lem}
\label{ptop^m}
If $f:X\ra Y$ is a $\mathbb Z/p$-homology equivalence, then it is a $\mathbb Z/p^N$-homology equivalence for all $N.$
\end{lem}
\begin{proof}
The long exact homology sequence associated to 

\begin{equation}
0\ra\mathbb Z\stackrel{\cdot p}\ra\mathbb Z\ra \mathbb Z/p\ra 0,
\end{equation}
implies that
$\tilde H_*(C_f;\mathbb Z/p)=0$ if and only if the map $\mathbb Z\stackrel{\cdot p}\ra\mathbb Z$
induces an isomorphism on $H_*(C_f;\mathbb Z).$ Then the map $\mathbb Z\stackrel{\cdot p^N}\ra\mathbb Z$
induces an isomorphism on $H_*(C_f;\mathbb Z)$ which, via the sequence

\begin{equation}
0\ra\mathbb Z\stackrel{\cdot p^N}\ra\mathbb Z\ra\mathbb Z/p^N\ra 0,
\end{equation}
implies that $\tilde H_*(C_f;\mathbb Z/p^N)=0.$
\end{proof}

\section{\label{Borel theorem}Borel's symmetry theorem for closed aspherical manifolds}

As a first application, we give a proof of the theorem of Borel (published by Conner and Raymond in \cite{connerraymond})
mentioned in the introduction.
The method in this simple situation is extended in the rest of this paper. 
\begin{thm}[3.2 in \cite{connerraymond}]
A closed asperical Riemannian manifold $(M,g)$ with centerless fundamental group $\pi_1M$
has no homotopically trivial isometries. 
\end{thm}
\begin{proof}
Isometries act on the fundamental group by moving loops around.
This gives a group homomorphism
\begin{equation}
\rho:\Isom(M,g)\ra\Out(\pi_1M).
\end{equation}  
Its kernel $K$ is the group of homotopically trivial isometries.
We will show that the kernel is trivial. 
By the Steenrod-Myers Theorem \cite{myerssteenrod}, the isometry group $\Isom(M,g)$ is a compact Lie group, 
so the closed subgroup $K$ is also a compact Lie group.
Since the center $Z(\pi_1M)$ is trivial, Proposition \ref{productlift}
implies that $K$ lifts to a subgroup $K<\Isom(\widetilde M,\widetilde g)$
which commutes with $\pi_1M.$ 
To show that the compact Lie group $K$ is trivial, it suffices to show that it has no elements of prime order $p$.
An element of order $p$ would define a $\mathbb Z/p$-action
which commutes with the $\pi_1M$ action. By Smith theory (Theorem \ref{coeffsmith}) we find that the
inclusion of the fixed point set of the $\mathbb Z/p$-action $F\hookrightarrow M$ is a $\mathbb Z/p$-homology equivalence.
In particular, it is an isomorphism
\begin{equation}
H_m(F;\mathbb Z/p)\cong H_m(M;\mathbb Z/p)
\end{equation}
on mod $p$ homology in dimension $m.$
By hypothesis, the $m$-dimensional manifold $M$ is compact. 
If $M$ is orientable, then it has a fundamental class $[M]\in H_m(M;\mathbb Z/p).$
This class corresponds to an $m$-dimensional class on the fixed point set $F$, so $F$ must be $m$-dimensional.
Since the fixed point set is a closed submanifold of $M,$ it must be all of $M$.
In other words, the $\mathbb Z/p$-action fixes everything, i.e. it is trivial.
If $M$ is not orientable, then we use Theorem \ref{coeffsmith} to get a $\mathbb Z/p$-homology
equivalence $F^{o}\hookrightarrow M^{o}$ in the orientable double cover $M^{o}\ra M$ and proceed as before
to conclude that the $\mathbb Z/p$-action is trivial. 
\end{proof}

\section{\label{intersections}Tubular neighborhoods and intersections of cycles}

In this section, we explain how to show that some manifolds are not bundles over lower dimensional submanifolds.
We need this for the following reason. The proof of Borel's theorem for compact manifolds reduced to showing
that there are no order $p$ homotopically trivial isometries. We used covering
space theory together with Smith theory to show that the inclusion of the fixed point set $F$ of such an isometry into
the manifold $M$ is a $\mathbb Z/p$-homology equivalence. When $M$ is compact, we saw that the fixed point
set is equal to all of $M$ using the fundamental class. 

Everything up to the last step works
the same way for a non-compact manifold $M$. To replace the last step, we proceed as follows.
The fixed point set of an isometry has a {\it tubular neighborhood} $NF$ which is a disk bundle over $F$. 
In detail, let
\begin{equation}
F\stackrel i\hookrightarrow NF\stackrel j\hookrightarrow M
\end{equation}
be the standard inclusions. Let $n$ be the codimension of the fixed point set $F$ inside the manifold $M$. 
Saying that $NF$ is a tubular neighborhood of $F$ means that we have a fibre bundle $\pi:NF\ra F,$ 
whose fibres are $n$-disks $\pi^{-1}(x)\cong D^n.$
The fixed point set $F$ is the zero section of this bundle, that is $\pi\circ i=\pi\mid_F=id_F.$ 

The tubular neighborhood $NF$ deformation retracts onto the fixed point set $F,$ 
so the inclusion $j:NF\hookrightarrow M$ is a $\mathbb Z/p$-homology equivalence.
To show that our isometry is trivial we need to show that the fixed point set is the entire manifold.
In other words we need to know that {\it$M$ is not $\mathbb Z/p$-homology
equivalent to a normal bundle of a lower-dimensional submanifold $F.$}  

Our goal now is to describe a homological condition which ensures this. We do this in 
subsection \ref{intersectioncondition}. To set things up, we recall
Poincare duality, intersections and the Thom 
isomorphism in the form we need (for non-compact manifolds).
A reference for the first two is \cite{hatcher}.

\subsection*{Notation}
For the remainder of this section, we adopt the following notations.
We have an orientable, $m$-dimensional manifold $M$
and an orientable submanifold $i:F\hookrightarrow M$ of lower dimension $m-n<m.$ 
The inclusion $i$ is closed, and the manifold $F$ has a tubular neighborhood $j:NF\hookrightarrow M.$
The tubular neighborhood is an $n$-disk bundle $\pi:NF\ra F$ and $i$ is the zero-section of this
bundle. Further, in the next three subsections, we will omit the coefficients from the notation 
for homology and cohomology groups.
The coefficients may be taken to be $\mathbb Q,\mathbb Z,$ or $\mathbb Z/p^N.$

\subsection{Poincare duality}
If $M$ is an orientable $m$-dimensional manifold, then we have a {\it Poincare duality isomorphism}
\begin{eqnarray}
PD:H_k(M)&\cong&H_c^{m-k}(M),\\
\ [a]&\mapsto&\eta_a,
\end{eqnarray}
from the homology of the manifold $M$ to the {\it compactly supported cohomology} in complementary dimension.
We denote the Poincare dual of a cycle $[a]$ by $\eta_a.$

\subsection{Cap products and intersections}
For any cohomology class $\omega\in H^n(M)$ and $k\geq 0,$ the {\it cap product with $\omega$} is a linear map

\begin{eqnarray}
H_{n+k}(M)&\stackrel{\omega\cap[-]}\longrightarrow &H_k(M),
\\
\ [b]&\mapsto&\omega\cap[b].\nonumber
\end{eqnarray}
Further, set $\omega\cap[b]=0$ when $[b]$ is a homology cycle of degree $<n.$ 

On the $n^{th}$ homology group, it is just defined by evaluating the cohomology class $\omega$ on the cycle $[b].$
On the higher degree homology groups ($k>0$) the cap product is defined on simplices by evaluating on only part of the
simplex. (See \cite{hatcher}.) A consequence of this is that if we can find representatives for $\omega$ and $[b]$
which have disjoint supports, then $\omega\cap[b]=0.$

For future reference we recall a formula relating cap and cup products. If $\alpha,\beta$ are
cohomology classes and $[c]$ is a homology cycle, then
\begin{equation}
\label{cupandcap}
(\alpha\cup\beta)\cap[c]=\alpha\cap(\beta\cap[c]).
\end{equation}

\vspace{0.5cm}
Let $[a]\in H_k(M)$ and $[b]\in H_{m-k}(M)$ be two cycles whose degrees add up to the dimension of $M.$
Call such a pair of cycles {\it complementary}. The 
{\it intersection product}\footnote{The left part of (\ref{intersectionproduct}) could be used to 
define the intersection product for more general pairs of cycles.
However, we will only use the notation $[a]\cap[b]$ when the cycles have complementary dimensions.} 
of two such cycles 

\begin{equation}
\label{intersectionproduct}
[a]\cap[b]:=\eta_a\cap[b]=\eta_a([b]),
\end{equation} 
is the Poincare dual of $[a]$ evaluated on the cycle $[b].$ 

On a connected manifold $M$, the intersection product has the following geometric interpretation.
Pick representatives for the cycles $[a]$ and $[b]$ which intersect transversally. 
The set theoretic intersection of these representatives is a 
finite collection of points. By keeping track of orientations, we assign each intersection a sign $\pm 1$. Then,
the intersection $[a]\cap[b]$ is the sum of these signs. 

\vspace{0.5cm}
{\it Cap products behave naturally with respect to open inclusions.}
Let $j: U\hookrightarrow V$ bet the inclusion of an open subset $U$, 
and let $\eta$ be an $n$-cocylce in $U$ which can be supported on a set 
$C\subset U$ that is closed in $V.$ Then, one can 
extend $\eta$ by zero to a cocycle $j_*\eta$ on $V.$ Let $[c]$ be a $k$-cycle. 
The definition of cap product on the level of simplices implies that 

\begin{equation}
\label{openinccap}
j_*(\eta\cap[c])=j_*\eta\cap j_*[c].
\end{equation}

If---in addition---$V$ is an orientable manifold, then intersection products on $U$ and $V$ make sense. 
For a pair of complementary cycles $[a]$ and $[b]$ in $U,$
we get

\begin{equation}
\label{openincint}
j_*([a]\cap[b])=j_*[a]\cap j_*[b].
\end{equation}
In other words, we can compute the intersection of the two cycles $[a]$ and $[b]$ either in the 
open submanifold $U$ or in the ambient manifold $V.$ 

Consequently, {\it if $j:U\hookrightarrow V$ is a
homology equivalence, then we can localize any intersection in $V$ to an intersection in $U,$}
because all the cycles in $V$ come from $U.$

\subsection{Thom isomorphism}
We use the notations from the beginning of this section. 
The tubular neighborhood is a bundle $\pi: NF\ra F.$ 
It has a {\it Thom class} 

\begin{equation}
Th(NF)\in H^n(NF,NF\setminus F).
\end{equation}

The fact that the Thom class lies in the relative cohomology group 
\newline 
$H^n(NF, NF\setminus F)$
means that its support can be localized to a smaller tubular neighborhood of the submanifold $F.$
Thus, the class is `compactly supported in the fibre direction.'

Multiplying a compactly supported class by the Thom class produces the {\it Thom isomorphism}
\begin{eqnarray}
H_c^{*}(F)&\stackrel{\mathrm{Thom}}\longrightarrow&H_c^{n+*}(NF),
\\
\eta&\mapsto& \pi^*\eta\cup Th(NF).\nonumber
\end{eqnarray}

The Thom class $Th(NF)$ is `compactly supported in the fibre direction' of the bundle $NF,$
while the cycle $\pi^*\eta$ is `compactly supported in the base direction'. Consequently, 
their product $Th(NF)\cup \pi^*\eta$ is compactly supported on the entire bundle $NF.$ 

The Thom isomorphism is related to Poincare duality by the following commutative diagram.

\begin{equation}
\begin{array}{ccc}
H^{d-n-*}_c(F)&\stackrel{\mathrm{Thom}}\longrightarrow&H^{d-*}_c(NF)
\\
\uparrow& &\uparrow
\\
H_*(F)&\stackrel{i_*}\longrightarrow &H_*(NF).
\end{array}
\end{equation}

Here, the vertical arrows are Poincare duality isomorphisms.

In other words, if $[a_0]\in H_k(F;\mathbb Z)$ is a homology cycle in the manifold $F$, then 
the Poincare dual of its image $i_*[a_0]$ inside $NF$ is given by the formula

\begin{equation}
\label{pdthom}
\eta_{i_*[a_0]}=\mathrm{Thom}(\eta_{a_0})=\pi^*\eta_{a_0}\cup Th(NF).
\end{equation}

Putting things together, we get the following lemma. 
It says that two complementary cycles in the tubular neighborhood $NF$ 
can be moved apart if one of them can be moved off the zero section $F$.

\begin{lem}
\label{thommove}
If $[a]$ and $[b]$ are complementary cycles and 
\begin{equation}
\label{moveoff}
Th(NF)\cap [b]=0
\end{equation} then the cycles can be moved apart. That is, 
\begin{equation}
[a]\cap[b]=0.
\end{equation}
\end{lem}
\begin{proof}
Since the bundle deformation retracts onto the zero section, there is a cycle $[a_0]\in H_{m-k}(F)$
such that $i_*[a_0]=[a].$ Thus, 
\begin{eqnarray*}
\ [a]\cap[b]&=&i_*[a_0]\cap [b],
\\
&=&\eta_{i_*[a_0]}\cap [b],\hspace{1.4cm}        \mbox{ def. of intersection}
\\
&=&(\pi^*\eta_{a_0}\cup Th(NF))\cap [b],\hspace{1.2cm}          \mbox{ by (\ref{pdthom})}
\\
&=&\pi^*\eta_{a_0}\cap(Th(NF)\cap[b]),\hspace{1.2cm}            \mbox{ by (\ref{cupandcap})}
\\
&=&0.\hspace{4.77cm}\mbox{by (\ref{moveoff})}
\end{eqnarray*}
\end{proof}


\subsection{\label{intersectioncondition}Showing that a manifold is not the normal bundle of a submanifold.}

We use the notations from the beginning of this section.
Consider the following condition on an $m$-dimensional, orientable manifold $M.$

\vspace{0.5cm}
\begin{itemize}
\item[($\pitchfork$)]
For any cohomology class $\omega\in H^n(M;\mathbb Z)$ of dimension $n>0$ 
there is some $k$ and complementary homology cycles 
$[a]\in H_{m-k}(M;\mathbb Z)$ and $[b]\in H_{k}(M,\mathbb Z)$
such that 
\begin{eqnarray}
\label{missF}
\omega\cap[b]&=&0,\\
\label{nonzeroint}
\ [a]\cap[b]&\not=&0.
\end{eqnarray}
\end{itemize}
\vspace{0.5cm}

Heuristically, this condition says that `$M$ has many intersections'.
We now show that such a manifold $M$ cannot look $\mathbb Z/p$-homologically like a 
bundle over a lower-dimensional submanifold. 

\begin{prp}
\label{notazpequivalence}
Let $M$ be an $m$-dimensional, orientable manifold and $F\hookrightarrow M$
an $(m-n)$-dimensional, orientable submanifold which is a closed subset of $M$. 
Let $j:NF\hookrightarrow M$ be a tubular neighborhood of this submanifold.
If $M$ satisfies condition $(\pitchfork)$ and $n>0,$ then the open inclusion $j$
is not a $\mathbb Z/p$-homology equivalence.
\end{prp}

\begin{proof}
The Thom class $Th(NF)$ can be supported on a smaller tubular neighborhood of $F,$
which means that it can be extended by zero to a class $j_* Th(NF)\in H^n(M;\mathbb Z)$ on the 
whole manifold $M.$
Applying 
($\pitchfork$) to $j_*Th(NF)$ gives a pair of cycles 
$[a]\in H_{m-k}(M;\mathbb Z)$ and $[b]\in H_k(M;\mathbb Z)$ such that 
\begin{eqnarray}
\label{cap}
j_*Th(NF)\cap[b]&=&0,
\\
\label{int}
\ [a]\cap[b]&\not=&0.
\end{eqnarray}

Pick a large enough power $N$ so that
$p^N$ does not divide the intersection number $[a]\cap[b].$ 
Then, equations (\ref{cap}) and (\ref{int})
hold in $\mathbb Z/p^N$-homology for the $\mathbb Z/p^N$-reductions of $j_*Th(NF), [a]$ and $[b].$

{\it Assume that $j$ is a $\mathbb Z/p$-homology equivalence.}
By Lemma \ref{ptop^m}, it is also a $\mathbb Z/p^N$-homology equivalence, 
so we find complementary cycles $[a_0],[b_0]\in H_*(NF;\mathbb Z/p^N)$
with $j_*[a_0]=[a]$ and $j_*[b_0]=[b].$ 

Plugging this in and using naturality of 
cap products and intersections under open inclusions (see (\ref{openinccap}) and (\ref{openincint})),
we get the equations
\begin{eqnarray*}
0&=&j_*Th(NF)\cap j_*[b_0]\\
&=&j_*(Th(NF)\cap[b_0]),\\
0&\not=&j_*[a_0]\cap j_*[b_0]\\
&=&j_*([a_0]\cap[b_0]),
\end{eqnarray*}
in $\mathbb Z/p^N$-homology. 
 
Since the map $j$ is a $\mathbb Z/p^N$ homology equivalence, we get
\begin{eqnarray*}
0&=&Th(NF)\cap [b_0],\\
0&\not=&[a_0]\cap [b_0].
\end{eqnarray*}
This contradicts Lemma \ref{thommove}.
Thus, $j$ cannot be a $\mathbb Z/p$-homology equivalence.
\end{proof} 



\section{\label{L2} Betti numbers of finite covers.}
In this section we produce a finite degree cover with
$M^i\ra M$ with `many intersections' (see Proposition \ref{nicecover}).
To do this, we first recall some facts about $L^2$-betti numbers. 
The examples that are of most interest to us are finite 
volume aspherical locally symmetric manifolds, and finite manifold covers of moduli space (i.e.
quotients of Teichm\"uller space by a finite index torsionfree subgroup of the mapping
class group). 
These examples are of the following sort. (See subsections \ref{locallysymmetric} and \ref{moduli} for more details and references.)
\begin{enumerate}
\item
The fundamental group $\Gamma:=\pi_1M$ is torsionfree and residually finite.
\item 
$M$ is the interior of a compact manifold with boundary $\partial M.$  
\item
The universal cover $\widetilde M$ is contractible and $\Gamma$-tame.
\item
The $L^2$-betti numbers of $\Gamma$ vanish except, possibly, in the middle dimension. 
The middle $L^2$-Betti number is equal to the absolute value of the Euler characteristic. 
\begin{equation}
b^k_{(2)}(\widetilde M;\Gamma)=\left\{\begin{array}{cc}
|\chi(M)| & \mbox{ if } k=m/2,
\\
0          & \mbox{ otherwise.}
\end{array}\right.
\end{equation}
\end{enumerate}
\begin{rem}
We do not know of any groups which satisfy the third condition but not the fourth one.
The question whether (3) implies (4) seems to be a sort of finite volume version of the Singer conjecture 
for closed aspherical manifolds.
\end{rem}
Let $\partial\widetilde M$ be the $\Gamma$-cover of the boundary.
Cheeger and Gromov showed 
that (2)+(3) implies that the boundary $\partial M$ is $L^2$-acyclic in the sense that 
\begin{equation}
\label{l2acyclic}
b_{(2)}^k(\partial \widetilde M;\Gamma)=0
\end{equation} 
for every $k.$
(This is Theorem $3.2$ of \cite{cheegergromov}
for residually finite $\Gamma$ and $1.2$ of \cite{cheegergromovbounds}
in general.)   
If $\Gamma$ is residually finite, then the following theorem of L\"uck lets 
one approximate the $L^2$-betti numbers by the (ordinary) betti numbers 
of a sequence of finite regular covers, normalized by the degree of the covers. 
\begin{thm}[L\"uck's approximation theorem \cite{luck1994}]
Let $\widetilde M$ be homotopy equivalent to a finite complex 
and $\widetilde M\ra M$ a regular cover
with residually finite group of covering translations $\Gamma$.
Pick regular finite covers $M^i\ra M$ of degree $d_i:=\mathrm{deg}(M^i\ra M)$
converging\footnote{This means $\pi_1\widetilde M=\cap_i\pi_1M^i$.} to $\widetilde M.$
Then,
\begin{equation}
b^k_{(2)}(\widetilde M;\Gamma)=\lim_{i\ra\infty}{b_k(M^i;\mathbb Q)\over d_i}.
\end{equation}
\end{thm}
The Euler characteristic is multiplicative in covers, $\chi(M^i)=d_i\cdot \chi(M),$ so 
L\"uck's approximation theorem implies 
\begin{equation}
\chi(M)=\chi^{(2)}(\widetilde M;\Gamma):=\sum_{k}(-1)^kb^k_{(2)}(\widetilde M;\Gamma).
\end{equation}

If $\chi(M)\not=0$ then there is a smallest non-zero $L^2$-betti number $b^k_{(2)}(\widetilde M;\Gamma).$
Since the fundamental group $\Gamma$ is infinite, the $d_i$ tend to infinity, so L\"uck's approximation
theorem together with (\ref{l2acyclic}) implies that there is a finite cover $M^i\ra M$ which satisfies

\begin{equation}
\label{bettiequation}
b_{k}(M^i;\mathbb Q)-b_{k-n}(M^i;\mathbb Q)>b_{k}(\partial M^i;\mathbb Q)
\end{equation}
for every $n>0.$ This is precisely what we need to prove Proposition \ref{nicecover}.
As a first step, we need the following lemma. It says that the cycles which intersect everything trivially
all come from the boundary.
\footnote{One direction of this is clear. If a cycle comes from the boundary then we can push
it towards the boundary and move it away from any other cycle.}

\begin{lem}
Let $(W,\partial W)$ be a $m$-dimensional compact manifold-with-boundary. Then, the image $I$ of the map
\begin{equation*}
H_{k}(\partial W;\mathbb Q)\ra H_{k}(W;\mathbb Q)
\end{equation*} 
consists of precisely those cycles $[a]\in H_{k}(W;\mathbb Q)$
for which $[a]\cap [b]=\eta_a\cap[b]=0$ for all $[b]\in H_{m-k}(W;\mathbb Q).$

The dimension of this space is $\dim_{\mathbb Q}I\leq b_{k}(\partial W;\mathbb Q).$
\end{lem}

\begin{proof}
Look at the commutative diagram
$$\begin{array}{cccccc}
\dots\ra &H^{m-k-1}(\partial W;\mathbb Q)&\stackrel{\delta}\ra &H^{m-k}(W,\partial W;\mathbb Q)&\stackrel{\tau}\ra &H^{m-k}(W;\mathbb Q)\ra\dots\\
         &\downarrow                   &    &\downarrow                     &     &\downarrow     \\
\dots\ra &H_{k}(\partial W;\mathbb Q)&\ra &H_k(W;\mathbb Q)&\ra &H_{k}(W,\partial W;\mathbb Q)\ra\dots.
\end{array}
$$
The horizontal lines are exact sequences in cohomology and homology, while the vertical lines are Poincare duality
isomorphisms. The subspace $I$ is Poincare dual to the kernel of $\tau,$ which consists of precisely those $\eta_a$ which
define trivial maps 
\begin{eqnarray*}
H_{k}(W;\mathbb Q)&\ra&\mathbb Q,\\
\ [b]&\mapsto&\eta_a\cap[b].
\end{eqnarray*}
This proves the first part of the lemma. The second part is clear.
\end{proof}

Next, we produce a finite cover with many intersections. 

\begin{prp}
\label{nicecover}
Let $M$ be an $m$-dimensional aspherical manifold with $\pi_1M$-tame universal cover. Suppose that
\begin{enumerate}
\item
The fundamental group $\pi_1M$ is residually finite.
\item 
$M$ is the interior of a compact manifold with boundary $\partial M.$  
\item
The Euler characteristic $\chi(M)$ is non-zero.
\end{enumerate}
Then, then there is a finite degree cover $M^i\ra M$
such that for any $n>0$ and $\omega\in H^n(M^i;\mathbb Z),$ there is $k$ 
and cycles $[a]\in H_{m-k}(M^i;\mathbb Z),[b]\in H_{k}(M^i;\mathbb Z),$
satisfying
\begin{eqnarray*}
\omega\cap[b]&=&0,\\
\ [a]\cap[b]&\not=&0.
\end{eqnarray*} 
\end{prp}

\begin{proof}
Let $M^i\ra M$ be a finite cover satisfying equation (\ref{bettiequation}).
Capping with the cohomology class $\omega$ gives a linear map
\begin{equation}
H_{k}(M^i;\mathbb Q)\stackrel{\omega\cap[-]}\longrightarrow H_{k-n}(M^i;\mathbb Q).
\end{equation}
Its kernel is a vector subspace of $H_{k}(M^i;\mathbb Q)$ of  dimension at least 
$b_{k}(M^i;\mathbb Q)-b_{k-n}(M^i;\mathbb Q)$ and by equation (\ref{bettiequation})
this number is greater than $b_{k}(\partial M^i;\mathbb Q).$ 
Thus, by the previous lemma, the kernel must contain a cycle $[b]_{\mathbb Q}$ on which the intersection form is non-zero. 
In other words, there are cycles of complementary dimensions $[a]_{\mathbb Q}$ and $[b]_{\mathbb Q}\in H_*(M^i;\mathbb Q)$ such that 
$\omega\cap[b]_{\mathbb Q}=0$ and $[a]_{\mathbb Q}\cap[b]_{\mathbb Q}\not=0$
in $H_*(M^i;\mathbb Q).$
Clearing denominators, we get cycles $[a]_{\mathbb Z},[b]_{\mathbb Z}\in H_*(M^i;\mathbb Z)$ so that

\begin{eqnarray}
\label{int1}
\omega\cap[b]_{\mathbb Z}&\mbox{ is }& \mbox{ torsion },\\
\label{int2}
\ [a]_{\mathbb Z}\cap[b]_{\mathbb Z}&\not=&0.
\end{eqnarray}
Multiplying $[b]_{\mathbb Z}$ by an appropriate positive integer $r$, we get a cycle $r[b]_{\mathbb Z}$ for which 
$\omega\cap r[b]_{\mathbb Z}=0$ while the intersection $[a]_{\mathbb Z}\cap r[b]_{\mathbb Z}$ is still nonzero.
\end{proof}

\subsection{\label{infinitetype} Manifolds of infinite type}
The assumption that $M$ is the interior of a compact manifold-with-boundary is not actually necessary.
Equations (0.4)-(0.13) in the paper \cite{cheegergromovbounds}
explain why, if $\widetilde M$ is $\pi_1M$-tameable, then the existence of a certain exhaustion $M=\cup_jM_j$ of $M$
by compact manifolds-with-boundary $M_j$ with controlled geometry of the boundary would imply
\begin{eqnarray*}
\lim_{j\ra\infty}b_{(2)}^k(\widetilde M_j;\Gamma)&=&b_{(2)}^k(\widetilde M;\Gamma),\\
\lim_{j\ra\infty}b_{(2)}^k(\partial\widetilde M_j;\Gamma)&=&0.
\end{eqnarray*}
Such an exhaustion is constructed in the later paper \cite{cheegergromovchoppings}.
Given this, one can establish Proposition \ref{nicecover} and all the remaining results
of this paper without assuming the existence of a boundary. In the statements of the theorems,
one should replace the Euler characteristic, which might not be defined when $M$
is not the interior of a manifold with boundary, by the $L^2$-Euler characteristic which
is always defined if $\widetilde M$ is $\pi_1M$-tameable.
Two arguments need to be modified.
\begin{itemize}
\item
For Proposition \ref{nicecover}, we use L\"uck's approximation theorem to find a cover $M_j^i\ra M_j$ such that 
\begin{equation*}
b_k(M_j^i;\mathbb Q)/b_0(M_j^i;\mathbb Q)-b_{k-n}(M_j^i;\mathbb Q)>b_k(\partial M_j^i;\mathbb Q)
\end{equation*}
for all $n>0$ and note that this implies some {\it connected} component $U_j^i$ of $M_j^i$
satisfies 
\begin{equation*}
b_k(U_j^i;\mathbb Q)-b_{k-n}(U_j^i;\mathbb Q)>b_k(\partial U_j^i;\mathbb Q)
\end{equation*}
for all $n>0.$ Then
proceed as in the proof of Proposition \ref{nicecover} to get the equations $\omega\cap[b]=0$ and $[a]\cap[b]\not=0$ on $U_j^i.$
The same equations hold on $M^i$ because $U_j^i\subset M_j^i$ is a connected open submanifold.
\item
For the compactness of the isometry group in the proof of Theorem \ref{maintheorem}
we show that for sufficiently large $j$, the subset $M_j$ cannot be moved off itself (because it
contains most of the $L^2$-cohomology) instead of $M_0.$  
\end{itemize}
Doing this has the advantage of avoiding the deep results about existence of boundaries
for locally symmetric manifolds and moduli spaces.


\section{\label{noactions}Proof of Theorem \ref{noactions2}: Eliminating periodic diffeomorphisms}
Let $\Gamma:=\pi_1M$ be the fundamental group of $M.$
We begin by proving the following part of the theorem.

\subsection*{Claim: If $\Gamma$ commutes with the compact Lie group $K$ then $K=1$}
To check that a compact Lie group is trivial, it suffices to check that it
has no elements of prime order $p$. Suppose we have such an element. In other words, 
we have a smooth $\mathbb Z/p$-action which commutes with the $\Gamma$-action. 
We apply Theorem \ref{coeffsmith} to show that in any cover $M'\ra M,$
the inclusion of the fixed point set $F'\hookrightarrow M'$ is a $\mathbb Z/p$-homology equivalence.
At this point, we want to pick $M'$ with `many intersections' (satisfying condition ($\pitchfork$))
and use Proposition \ref{notazpequivalence} to show that the fixed point set $F'$ is all of $M'$. 
However, to apply this proposition, we first need to make sure that the manifold $M'$ {\it and} the fixed point set $F'$
are orientable.
If $M$ is not orientable, let $M^o\ra M$ be the orientation double cover. (Set $M^o=M$ is $M$ is orientable.)
\begin{itemize}
\item
If $p$ is odd, then the $\mathbb Z/p$-action is orientation-preserving, so that the fixed point set $F^o$ is orientable.
\item
If $p=2,$ then we know that the inclusion of the fixed point set $F^o\hookrightarrow M^o$ is a $\mathbb Z/2$-homology
equivalence. Thus, it induces an isomorphism 
$H^1(M^o;\mathbb Z/2)\cong H^1(F^o;\mathbb Z/2).$
If $F^o$ is not orientable, then its orientation double cover is classified by a map $F^o\ra K(\mathbb Z/2,1)$. The isomorphism
above shows that this map factors through $F^o\hookrightarrow M^o\ra K(\mathbb Z/2,1).$ Thus, taking a further
double cover of $M^o$ if necessary, we may assume that $F^o$ is orientable. 
\end{itemize}
Now, the manifold $M^o$ is satisfies the conditions of Lemma \ref{nicecover} 
so it has a further finite degree cover $M'\ra M^o\ra M$ satisfying 
($\pitchfork$). 
Finally, the fixed point set $F'$ and the cover $M'$ are orientable, because they cover $F^o$ and $M^o$, respectively.
At this point, we can apply Proposition \ref{notazpequivalence}. Since $F'\hookrightarrow M'$ is a $\mathbb Z/p$-homology
equivalence, and $M'$ has `many intersections' we conclude that the fixed point set $F'$ is all of $M'.$ 
This shows that the $\mathbb Z/p$-action is trivial.
Consequently, the group $K$ is trivial. This proves the claim. 

\subsection*{No homotopically trivial periodic diffeomorphisms}
Since $\chi(M)\not=0$ the manifold $M$ has centerless fundamental group (7.2 in \cite{luck}),
so Proposition \ref{productlift} shows that a homotopically trivial periodic diffeomorphism $f$ lifts
to a periodic diffeomorphism $\widetilde f$ of the universal cover which commutes with the fundamental group. The above
claim shows that this lift $\widetilde f$ is the identity, so the diffeomorphism $f$ is also the identity map. 

\subsection*{If $\Gamma$ normalizes a compact Lie group $K$, then $K=1$} 
We do this by successively eliminating characteristic compact subgroups. This is similar to what is 
done in \cite{farbweinbergerroyden}. 
Let $K_0<K$ be the identity component of the compact Lie group. We have the exact sequence
\begin{equation}
1\ra K_0^{sol}\ra K_0\ra K_0^{ss}\ra 1,
\end{equation}
where $K_0^{sol}$ is the maximal normal connected solvable subgroup and $K_0^{ss}$ is semi-simple. 
Since $K_0^{sol}$ is a compact connected solvable group, it is a torus $(S^1)^{n}.$ 
The finite subgroup $(\mathbb Z/p)^n<(S^1)^n$ is characteristic (it consists of the elements of order $1$ and $p$ in the torus),
$K_0^{sol}$ is characteristic in $K_0,$ the group $K_0$ is (topologically) characteristic in $K$ and $\Gamma$ acts on $K$ by conjugation, 
so $\Gamma$ also acts on the group $(\mathbb Z/p)^n$ by conjugation. Since this last group is finite, some finite index subgroup $\Gamma'<\Gamma$
acts trivially on it, i.e. $\Gamma'$ commutes with $(\mathbb Z/p)^n.$ By the claim above,
we have $n=0$ so that $K_0$ is semisimple. Thus, the center $Z(K_0)$ is finite so some finite index subgroup of $\Gamma$
commutes with it and by the claim we find that the center is trivial. Let $\Gamma'<\Gamma$ be the kernel of the conjugation homomorphism $\Gamma\ra\Out(K_0).$  Since $K_0$ is semisimple,
its outer automorphism group $\Out(K_0)$ is finite so the kernel is a finite index subgroup. By Proposition \ref{productlift},
the group generated by $K_0$ and $\Gamma'$ splits as a product $K_0\times\Gamma'$ and by the claim above
we conclude that $K_0$ is trivial. We've shown that $K$ is a finite group. Thus, there is a finite index subgroup $\Gamma'<\Gamma$
which commutes with $K$ and by the claim $K$ is trivial.

\section{\label{recognition}Isometry groups of universal covers}
In this section, we prove Theorem \ref{recognizinglocsym}.
In fact, we prove the following generalization of that theorem, which allows
some torsion in the covering group $\Gamma$ and realizes the isometry group 
of the universal cover as a subgroup of an
algebraic gadget called the abstract commensurator of $\Gamma.$
\begin{defn}
The {\it abstract commensurator} $\Comm(\Gamma)$ is the group of equivalence classes of 
isomorphisms $\psi:\Gamma_1\ra\Gamma_2$ between finite index subgroups $\Gamma_1, \Gamma_2$ of $\Gamma$,
with two such isomorphisms $\psi:\Gamma_1\ra\Gamma_2$ and $\psi':\Gamma_1'\ra\Gamma_2'$
being equivalent if they agree on some further finite index subgroup of $\Gamma.$ 
\end{defn} 

\begin{thm}
\label{recognize2}
Let $\widetilde M$ be a contractible $\Gamma$-tame manifold. 
Suppose the action of $\Gamma$ on $\widetilde M$ is effective and there is a finite index torsionfree subgroup $\Gamma'<\Gamma.$ 
Suppose, in addition, that the following three conditions hold.
\begin{enumerate}
\item
$\Gamma$ is residually finite and irreducible\footnote{That is, if $A\times B$ is a finite index subgroup of 
$\Gamma$ then either $A$ or $B$ is finite.}. 
\item
$\widetilde M/\Gamma'$ is the interior of a compact manifold with boundary.
\item
The Euler characteristic $\chi(\widetilde M/\Gamma')$ is non-zero.
\end{enumerate}
If $\widetilde g$ is a complete Finsler metric on $\widetilde M$ which is $\Gamma$-invariant and has finite
$\Gamma$-covolume
then either
\begin{itemize}
\item
$(\widetilde M,\widetilde g)$ is isometric to a symmetric space, 
or
\item
The isometry group $\Isom(\widetilde M,\widetilde g)$ is discrete and $\Gamma$ is a finite index subgroup.
Moreover, there exist inclusions
\begin{equation}
\label{commensurator}
\Gamma<\Isom(\widetilde M,\widetilde g)<\Comm(\Gamma).
\end{equation}
\end{itemize}
\end{thm}
The conditions of the theorem imply that the group $\Gamma$ is centerless (see the remark below.)
The substance of equation (\ref{commensurator}) is that the isometry group 
injects into the abstract commensurator. This is a group theoretic analogue
of the geometric statement that there are no homotopically 
trivial isometries in any finite cover of $\widetilde M/\Gamma'$.
It means that information about the commensurator can give 
quantitative bounds on the size of the isometry group.
In general, the abstract commensurator is difficult to compute, 
but much is known about it for lattices in symmetric spaces,
and it is completely determined for (finite index subgroups of)
mapping class groups. 

The proof below uses methods developed by Farb and Weinberger in \cite{farbweinbergerroyden}.
The main new ingredient is Theorem \ref{noactions2}, which allows us to eliminate 
compact Lie groups commuting with the covering group. Previously, one only
had a way of eliminating {\it connected} compact Lie groups 
which was not enough to get a quantitative bound on the size of the isometry group.
Moreover, the previous method relied on Despotovic's result \cite{despotovic} about the action dimension
of the mapping class group, while our method does not. We also replace a couple arguments
specific to the mapping class group by general $L^2$-cohomology arguments. 
\begin{itemize}
\item
If $\chi^{(2)}(\Gamma)\not=0$ then $\Gamma$ has no infinite normal amenable subgroups (7.2 (1), (2) in \cite{luck}.)
\item
If $\chi^{(2)}(\Gamma)\not=0$, then any two contractible $\Gamma$-tame manifolds $\widetilde M$ and $\widetilde N$
have the same dimension (5.2 in \cite{cheegergromovgroups}). 
\end{itemize}
\begin{rem}
Let $\Gamma$ be the group in Theorem \ref{recognizinglocsym}.
Since it has non-zero $L^2$-Euler characteristic, its
center $Z(\Gamma)$ is finite. Applying Theorem \ref{noactions2} to the quotient of 
$\widetilde M$ by a finite index torsionfree subgroup $\Gamma'<\Gamma$ we find
that the center acts trivially on $\widetilde M.$ Since we've assumed that $\Gamma$
acts effectively, we conclude that $Z(\Gamma)=1.$
\end{rem}
\subsection*{Proof of Theorem \ref{recognize2}}
Let $I:=\Isom(\widetilde M,\widetilde g)$ be the isometry group and $I_0$ its identity component.
The group $I$ is a Lie group acting properly and smoothly on $\widetilde M.$ (\cite{myerssteenrod}
in the Riemannian case and \cite{denghou} in the general Finsler case.)
Let $\Gamma_0=\Gamma\cap I_0.$ 
Then, we have an exact sequence 
\begin{equation}
1\ra I_0^{sol}\ra I_0\ra I_0^{ss}\ra 1
\end{equation} 
where $I_0^{sol}$ is the maximal connected solvable normal subgroup of $I_0$ and $I_0^{ss}$ is semi-simple.
The group $\Gamma_0$ has \cite{raghunathan} a unique maximal normal solvable subgroup $\Gamma_0^{sol}.$
This subgroup is characteristic in $\Gamma_0$, so it is normal in $\Gamma.$ Since $\Gamma_0^{sol}$ is 
solvable, it is amenable, so $\chi^{(2)}(\Gamma)\not=0$ implies this group is finite.
Farb and Weinberger \cite{farbweinbergerroyden} show that $\Gamma_0$ is a lattice in $I_0$. (This
only uses the assumption that $\widetilde g$ has finite $\Gamma$-covolume and
no particular properties of the mapping class group situation.)
This lets them use Prasad's Lemma $6$ in \cite{prasad} to conclude that $I_0^{sol}$ 
is compact and the center $Z(I_0^{ss})$ is finite. Since $I_0^{sol}$ is characteristic,
it is normalized by $\Gamma$. Since it is compact, Theorem \ref{noactions2} implies it is trivial. Thus, $I_0=I_0^{ss}$.
Now, $\Gamma$ acts by conjugation on the finite group $Z(I_0^{ss})$ so again Theorem \ref{noactions2}
shows $Z(I_0^{ss})=1.$ Thus, $I_0$ is semisimple with trivial center.
Next, we note that it cannot have any compact factors: the product $K$ of the compact factors
is a characteristic subgroup, so it is normalized by $\Gamma$ and by Theorem
\ref{noactions2} it must be trivial.
Thus, $I_0$ is semisimple with trivial center and no compact factors.
Now, look at the extension
\begin{equation}
1\ra I_0\ra<I_0,\Gamma>\ra\Gamma/\Gamma_0\ra 1.
\end{equation}
Let $\Gamma'<\Gamma$ be the kernel of the conjugation homomorphism $\Gamma\ra\Out(I_0).$ 
It is a finite index subgroup of $\Gamma$, since $I_0$ is semisimple with trivial center. By
Proposition \ref{productlift} the extension $<I_0,\Gamma'>$ splits as
$I_0\times (\Gamma'/\Gamma'_0),$ where $\Gamma'_0=\Gamma'\cap I_0.$ Consequently, 
$\Gamma'=\Gamma'_0\times (\Gamma'/\Gamma'_0)$. Since $\Gamma'$ is irreducible, we must have either
that $\Gamma'_0$ is finite or that $\Gamma'_0$ is a finite index subgroup of $\Gamma'.$ Now we look at these two possibilities.
\subsection*{Case 1: $\Gamma'_0$ is a finite index subgroup of $\Gamma'$.}
Recall that $\Gamma_0$ is a lattice in the Lie group $I_0.$
Since $\Gamma'_0$ is a finite index subgroup of $\Gamma_0,$ it is also a lattice in $I_0.$
Let $K<I_0$ be a maximal compact subgroup. The quotient $I_0/K$ is a symmetric space
with no compact or Euclidean factors, so it is contractible and $\Gamma'_0$-tame (see subsection \ref{locallysymmetric}).
The manifold $\widetilde M$ is also contractible and $\Gamma$-tame. Since $\Gamma'_0<\Gamma$ is a finite index subgroup,
$\widetilde M$ is also $\Gamma'_0$-tame. Since $\chi^{(2)}(\Gamma'_0)\not=0,$ $I_0/K$ and $\widetilde M$ have 
the same dimension. (See Corollary 5.2 in \cite{cheegergromovgroups}).
On the other hand, an orbit of $I_0$ acting on $\widetilde M$ has the form $I_0\cdot x=I_0/K_x$ for some compact subgroup $K_x<I_0.$
We have $\dim\widetilde M=\dim I_0/K\leq\dim I_0/K_x\leq\dim\widetilde M$ so the inequalities are actually
equalities. In other words, the isometry group $I_0$ is acting transitively on $\widetilde M$, i.e. 
$(\widetilde M,\widetilde g)$ is isometric to a symmetric space.  
\subsection*{Case 2: $\Gamma'_0$ is finite.}
If $\Gamma'_0$ is finite, so is $\Gamma_0.$ Since $\Gamma_0$ is a lattice in $I_0$ and $I_0$ has no compact factors, 
this implies that $I_0$ is trivial, so that the isometry group $I$ is discrete. 
Since $\Gamma$ is a lattice in $I$, we conclude that it is a finite index
subgroup of $I.$ Conjugation by elements of $I$ preserves some finite index subgroup $\Gamma'<\Gamma$,
so it defines a commensuration of $\Gamma.$ This gives a homomorphism
\begin{equation}
I\ra\Comm(\Gamma).
\end{equation} 
Let $f$ be in the kernel of this homomorphism. It acts trivially by conjugation on some 
subgroup $\Gamma'<\Gamma$, so it commutes with that subgroup. Since $\Gamma'$ 
has trivial center ($\chi^{(2)}(\Gamma')\not=0$ implies the center is 
finite, and a finite group commuting with $\Gamma'$ must be trivial by
Theorem \ref{noactions2}), we have $<f>\cap\Gamma'=1$ so that the group $<f,\Gamma'>$
is a product $<f>\times\Gamma'$ and the element $f$ has order at most $|I/\Gamma'|<\infty.$
By Theorem \ref{noactions2}, $f$ is trivial, so that $I$ is a subgroup of $\Comm(\Gamma).$

This completes the proof of Theorem \ref{recognize2} and of its special case Theorem \ref{recognizinglocsym}.
\begin{rem}
For future use, we also note that if $\widetilde M/\Gamma$ is a finite branched cover of an orbifold $Y$,
then the orbifold fundamental group $\pi_1^{orb}Y$ acts properly discontinuously on $\widetilde M$ 
and contains $\Gamma$ as a finite index subgroup. The argument above shows that there are inclusions
\begin{equation}
\label{branchedcover}
\Gamma<\pi_1^{orb}Y<\Comm(\Gamma).
\end{equation}
This can sometimes be used to bound the 
degree of the branched cover $\widetilde M/\Gamma\ra Y.$
\end{rem}
\section{\label{royden} Moduli spaces}
\subsection{\label{moduli}Background on moduli spaces}
A common reference for most of the properties recalled here is the survey \cite{ivanov}.
Let $S=S_{g,n}$ be the genus $g$ surface with $n$ punctures and $\Mod(S)$ its mapping class group. 
The Teichm\"uller space $\mathrm{Teich}(S_{g,n})$ has real dimension $6g-6+2n.$
The mapping class group is residually finite \cite{grossman} and virtually torsionfree.
McMullen's K\"ahler hyperbolic metric is a complete Riemannian bounded geometry metric
which is $\Mod(S)$-invariant and has finite $\Mod(S)$-covolume \cite{mcmullen}.
Consequently, Teichm\"uller space is $\Mod(S)$-tame. Since the McMullen metric is K\"ahler hyperbolic,
the Euler characteristic of the mapping class group is non-zero. 
Teichm\"uller space has   
a partial compactification $\overline T$ (analogous to the Borel-Serre
partial compactification of arithmetic symmetric spaces) for which $\overline T/\Gamma$ is a compact manifold-with-boundary
for any finite index torsionfree subgroup $\Gamma<\Mod(S)$.\cite{harvey,ivanov2}
The mapping class group is irreducible (c.f. page $8$ of \cite{farbweinbergerroyden}).
\begin{rem}
Below, we will only look at the moduli spaces of real dimension at least four.
(The zero dimensional moduli spaces are finite, and the two dimensional moduli spaces are hyperbolic.)
\end{rem}

\subsection{Proof of Theorem \ref{generalroyden2}}
The above paragraph implies that Theorem \ref{recognize2}
applies to finite index subgroups $\Gamma<\Mod^{\pm}(S)$ 
of the extended mapping class group (modulo the center in the cases $S_{1,2}$ and $S_{2,0}$) 
acting on Teichm\"uller space. Since such a $\Gamma$ is not a lattice in any symmetric space
(1.2a in \cite{farbweinbergerroyden}), we conclude that the isometry group $I:=\Isom(\mathrm{Teich}(S),\widetilde g)$
is discrete and contains $\Gamma$ as a finite index subgroup.

Let $i:\Gamma\ra\Comm(\Gamma)$
be the action of $\Gamma$ on itself by conjugation. 
If $S$ is not a sphere with $\leq 4$ punctures and not a torus with $\leq 2$ punctures,  then
$\Comm(\Gamma)=i(\Mod^{\pm}(S))$. (8.5A and 9.2B in \cite{ivanov})
If $S_{1,2}$ is the torus with two punctures, then $i(\mathrm{Mod}^{\pm}(S_{1,2}))$
is an index five subgroup of the commensurator $\Comm(\mathrm{Mod}^{\pm}(S_{1,2}))\cong\mathrm{Mod}^{\pm}(S_{0,5})$ 
(see Proposition $8$ in \cite{margalit}.) 
Plugging this information into equation (\ref{commensurator}) 
shows that 
\begin{equation}
\label{bound}
|I/\Gamma|\leq\left\{\begin{array}{cc}|\Mod^{\pm}(S)/\Gamma|&\mbox{ if } S\not=S_{1,2},\\
5|\Mod^{\pm}(S)/\Gamma|&\mbox{ if } S=S_{1,2}.\end{array}\right.
\end{equation}
This proves Theorem \ref{generalroyden2} when $\Gamma=\Mod^{\pm}(S),$
and also gives bounds for the isometry groups of metrics
that are only invariant under a finite index subgroup of $\Mod^{\pm}(S).$

\subsection{Isometries of the Teichm\"uller metric}
Royden proved that the isometry group of the Teichm\"uller
metric on the Teichm\"uller space of a closed genus $g\geq 2$
surface is the extended mapping class group\cite{royden}
and Earle-Kra\cite{earlekra} extended this result to genus
$g$ surfaces with $n$ punctures whenever $6g-6+2n\geq 4$ and $(g,n)\not=(1,2).$
In the exceptional case $S_{1,2}$, they showed that the Teichm\"uller spaces
of $S_{1,2}$ and $S_{0,5}$ are isometric in the Teichm\"uller metric,
and consequently the isometry group of $\mathrm{Teich}(S_{1,2})$
in that metric is the extended mapping class group of $S_{1,5}$.
Thus, the inequality (\ref{bound}) is sharp in all cases.


\subsection{Minimal orbifold} The information about the commensurator
of the mapping class group described above, together with equation (\ref{branchedcover}), 
shows that $\mathcal M_{g,n}^{\pm}$ is a minimal orbifold unless $(g,n)=(1,2).$
\section{\label{locsym} Locally symmetric spaces}
\subsection{\label{locallysymmetric}Background on locally symmetric spaces}
If $(M,g)$ is a closed Riemannian manifold, then its sectional curvatures are bounded,
and its injectivity radius is positive, so it has bounded geometry. Thus, its universal cover $(\widetilde M,\widetilde g)$
also has bounded geometry. 
Any symmetric space $(G/K,\widetilde g_{sym})$ 
with its standard symmetric metric has bounded geometry because it isometrically covers a 
compact locally symmetric space. Thus, it is $\Gamma$-tame for every lattice $\Gamma<G.$
Any finite volume locally symmetric space $M=\Gamma\setminus G/K$ is the interior of a compact
manifold-with-boundary \cite{borelserre,borelji,pettetsouto}. 
When $M$ is an {\it aspherical} locally symmetric space, 
one knows that the middle-dimensional 
$L^2$-Betti number is $|\chi(M)|$, and all other $L^2$-Betti
numbers vanish (see \cite{olbrich}; he cites \cite{borel} as the original reference).

\subsection{Proof of Theorem \ref{symmetric}}
The above paragraph shows that Theorem \ref{recognize2}
applies to an irreducible $\chi\not=0$ lattice $\Gamma$ in the semisimple
Lie group $G:=\Isom(\widetilde M,\widetilde h_{sym}).$
If $\widetilde g$ is a constant multiple of the symmetric metric, then there
is nothing to prove, so assume that it is not. 
In this case $I:=\Isom(\widetilde M,\widetilde g)$
is discrete, and contains $\Gamma$ as a finite index subgroup.

\vspace{0.2cm}
\noindent{\bf Claim:} There is an effective, properly discontinuous action of the group $I$
which preserves the symmetric metric $\widetilde h_{sym}$ on $\widetilde M.$ 
\vspace{0.2cm}

Given this claim, we get 
\begin{equation*}
|I/\Gamma|={\vol(\widetilde M/\Gamma,\widetilde h_{sym})\over\vol(\widetilde M/I,\widetilde h_{sym})}.
\end{equation*}
Kazhdan and Margulis have shown that the volume of the smallest locally
symmetric orbifold covered by a symmetric space with no compact of Euclidean factors $(\widetilde M,\widetilde h_{sym})$
is bounded below by a positive constant $\varepsilon(\widetilde h_{sym})$ 
that depends only on the symmetric space (c.f. XI.11.9 in \cite{raghunathan}).
This proves Theorem \ref{symmetric}, given the claim.
\begin{rem}
By contrast, note that the real line covers circles of arbitrarily small volume. 
\end{rem}

Now we prove the claim. Since $\Gamma<I$ is a finite index subgroup, there is a further finite index subgroup
$\Gamma'<\Gamma$ which is normalized by $I$. Moreover, since $\Gamma$ is virtually torsionfree,
we may take $\Gamma'$ to be torsionfree.
Now, conjugation by $I$ gives a group homomorphism
\begin{equation}
\rho:I/\Gamma'\ra\Out(\Gamma')
\end{equation}
and it is injective because $\widetilde M/\Gamma'$ has no homotopically trivial periodic diffeomorphisms.
The image $\rho(I/\Gamma')$ can be represented by isometries of a locally symmetric metric
(by Kerckhoff's proof of Nielsen realization\cite{kerckhoff} when $\widetilde M/\Gamma'$ is a surface, 
and by Margulis-Mostow-Prasad rigidity if $\widetilde M/\Gamma'$ is an irreducible locally symmetric
space of dimension $>2$.)
Thus, the group $L$ of lifts
\begin{equation*}
1\ra\Gamma'\ra L\ra\rho(I/\Gamma')\ra 1
\end{equation*} 
can be represented by isometries of the symmetric metric. This extension is determined
by $\rho$ since $\Gamma'$ has trivial center (see section \ref{lifts})
so the group of lifts $L$ is isomorphic to $I$. This proves the claim.

\subsection{Most symmetry}
By Theorem \ref{maintheorem} the isometry
group of any complete Riemannian metric $g$ on a $\chi\not=0$
locally symmetric space $M$ is isomorphic to a finite subgroup of $\Out(\pi_1M).$
If $M$ is irreducible then, by Nielsen realization/Mostow rigidity, 
the isometry group is isomorphic to a subgroup of the isometry group of a 
locally symmetric metric $\Isom(M,g)<\Isom(M,h_{sym}).$
Thus, in a sense, the locally symmetric metric has the most symmetry. 
\section{\label{BorelStar}Infinite volume metrics and isometries of the base: Proof of Theorem \ref{maintheorem}}
The main point is to show that the group $\Isom(M,g)$ is compact, even if the metric $g$ has infinite volume\footnote{For a finite volume complete metric $g$, this is standard.}.
Then we will show that there are no homotopically trivial isometries.
\subsection*{Compactness of the isometry group}
Let $M_0\subset M$ be the complement of an open collar neighborhood of the boundary.
We first show that this subset cannot be moved off itself by a homeomorphism because it
`contains all of the $L^2$-cohomology'.
If a homeomorphism $\phi$ moves $M_0$ completely off itself ($\phi(M_0)\cap M_0=\emptyset$),
then the induced isomorphism $\phi^*$ on $L^2$-homology factors through $M\setminus M_0\sim \partial M.$
That is, we have maps 
\begin{equation*}
H^*_{(2)}(\widetilde M;\Gamma)\ra H^*_{(2)}(\partial\widetilde M;\Gamma)\ra H^*_{(2)}(\widetilde M;\Gamma)
\end{equation*}
and the composition is the isomorphism $\phi^*$. This cannot happen because $\partial M$ is $L^2$-acyclic while $M$ is not. 

We've shown that the subset $M_0\subset M$ cannot be moved off itself.
Thus, a point $x\in M_0$ gets moved at most a distance $2\cdot\mathrm{diameter}(M_0)=:2D$
by any isometry of $(M,g).$ The isometry group $\Isom(M,g)$ is a Lie group acting smoothly on $M$,
and the action
\begin{eqnarray*}
A:\Isom(M,g)&\ra& M\times M\\
(\phi,x)&\mapsto&(\phi(x),x)
\end{eqnarray*}
is a proper map (\cite{myerssteenrod} in the
Riemannian case and \cite{denghou} in the general Finsler case). Thus $A^{-1}(\overline{B_x(2D)}\times\{x\})$ is compact. 
It contains the closed subset $\Isom(M,g)\times\{x\}$ (because the orbit $\Isom(M,g)\cdot x$ is contained
in the closed ball $\overline{B_x(2D)}$) which shows that $\Isom(M,g)$ is compact. 

\subsection*{No homotopically trivial isometries}
Let $K:=\ker(\Isom(M,g)\ra\Out(\pi_1M))$ be the group of homotopically trivial isometries.
It is a closed subgroup of the isometry group so it is also a compact Lie group. 
Since $M$ has no periodic homotopically trivial isometries (Theorem \ref{noactions2}),
this group is trivial.
Thus, $\Isom(M,g)<\Out(\pi_1M)$
is a compact subgroup of a discrete group. Consequently, it is finite. 

\section{\label{infinitevolume}Infinite volume metrics and isometries of the universal cover}
In this section, we give an example which shows that 
Theorem \ref{recognizinglocsym} does not extend to metrics $g$ of 
infinite volume.

Let $M$ and $N$ be non-compact aspherical manifolds whose universal cover is 
diffeomorphic to $\mathbb R^m.$ Pick a proper embedded rays $r_1:[0,\infty)\ra M$
and $r_2:[0,\infty)\ra N$, cut out tubular neighborhoods of these and glue 
to get a new manifold $M\#_rN$ which is a connect sum along the rays. 
This manifold is homotopy equivalent to the wedge $M\vee N,$ so it is aspherical,
and its fundamental group is $\pi_1M\star\pi_1N.$ The quotient $M':=(\widetilde{M\#_rN})/\pi_1M$
is homotopy equivalent to $M$ and, in fact, one can check that it is diffeomorphic to $M.$ 
Thus, $M$ is an infinite cover of the manifold $M\#_rN.$ Now, pick a complete Riemannian metric $\overline g$
on $M\#_rN.$ It lifts to an infinite volume complete Riemannian metric $g$ on $M,$ which further 
lifts to a metric $\widetilde g$ on the universal cover whose isometry group $\Isom(\widetilde M,\widetilde g)$
contains the free product $\pi_1M\star\pi_1N.$ 


\bibliography{locsym}
\bibliographystyle{alpha}

\end{document}